\documentclass[11pt,a4paper]{article}
\usepackage[utf8]{inputenc}
\usepackage[english]{babel}
\usepackage{amsmath}
\usepackage{amsfonts}
\usepackage{amssymb}
\usepackage{makeidx}
\usepackage{graphicx}
\usepackage{lmodern}
\usepackage[left=2cm,right=2cm,top=2cm,bottom=2cm]{geometry}

\usepackage{dsfont}
\usepackage{mathtools}
\usepackage{tikz}
\usepackage{relsize}
\usepackage{amsmath,environ}
\usepackage[symbol]{footmisc}

\usepackage{calrsfs}
\DeclareMathAlphabet{\pazocal}{OMS}{zplm}{m}{n}

\newcommand{\R}{\mathbb{R}} 
\newcommand{\C}{\mathbb{C}}

\newcommand{\Sbb}{\mathbb{S}}

\newcommand{\Hpazo}{\pazocal{H}}

\newcommand{\Dcal}{\mathcal{D}}
\newcommand{\Ecal}{\mathcal{E}}

\newcommand{\Lcal}{\mathcal{L}}

\newcommand{\Ccal}{\mathcal{C}}
\newcommand{\Ocal}{\mathcal{O}}
\newcommand{\Scal}{\mathcal{S}}

\newcommand{\Lip}{\textnormal{Lip}}
\newcommand{\Span}{\textnormal{Span}}

\newcommand{\Self}{\textnormal{Self}}
\newcommand{\ASelf}{\textnormal{A-Self}}
\newcommand{\Vecf}{\textnormal{Vec}}
\newcommand{\E}{\textnormal{E}}

\newcommand{\Conj}{\textnormal{Conj}}
\newcommand{\g}{\mathbf{g}}

\newcommand{\INTSeg}[4]{\int_{#3}^{#4} #1 \textnormal{d} #2}

\newtheorem{lem}{Lemma}
\newtheorem{Def}{Definition}

\newtheorem{thm}{Theorem}
\newtheorem{prop}{Proposition}
\newtheorem{cor}{Corollary}

\author{Beno\^it Bonnet\footnote{Aix Marseille Universit\'e, CNRS, ENSAM, Universit\'e de Toulon, LIS, Marseille, France.} , Jean-Paul Gauthier$^*$, Francesco Rossi\footnote{Dipartimento di Matematica "Tullio Levi-Civita" Universit\`a degli Studi di Padova, Padova, Italy.}} 

\title{Generic Singularities of the 3D-Contact Sub-Riemannian Conjugate Locus}

\begin{document}
\maketitle

\begin{abstract}
In this paper, we extend and complete the classification of the generic singularities of the 3D-contact sub-Riemmanian conjugate locus in a neighbourhood of the origin. 
\end{abstract}

\section{Introduction}
\label{section:Introduction}

One of the major discrepancies between Riemannian and sub-Riemannian geometry is the fact that the sub-Riemannian distance generically fails to be smooth in a neighbourhood of its base point. This uncanny property translates into the presence of a wide amount of singularities along the so-called \textit{caustic} surfaces generated by the critical points of the corresponding exponential map. The analysis of the least degenerate generic behaviour of the 3D-contact sub-Riemannian caustic was carried out for the first time in the seminal paper \cite{A1}. The further degenerate generic situations were then studied independently in \cite{A2,CGK}, and the full classification of the generic singularities of the \textit{semi-caustics} - i.e. the intersections of the caustic with adequate half-spaces - was completed in \cite{ACG}. In the present note, we accomplish this research by providing a complete classification of the generic singularities of the full caustic. Let it be noted that the question of \textit{stability} for these singularities is very delicate, and was only studied subsequently in \cite{ACGZ1}. Furthermore, the analysis of the less degenerate configurations of the contact sub-Riemannian caustics has been recently extended to the case of dimension greater or equal to 5 in \cite{S}. 

It was shown in \cite{ACG,CGK} that for the caustic, the generic degenerate situations appear when one of the fundamental geometric invariants of the sub-Riemannian distribution vanishes. Starting from the generic non-degenerate situation in which the semi-caustics exhibit 4 plea lines (see Figure 1-left), one can observe the formation of a \textit{swallow tail} (see Figure 1-center) as the geometric invariant gets closer to zero. When this invariant vanishes, the swallow tail folds itself and becomes degenerate, and the semi-caustics become distinct closed surfaces exhibiting 6 plea lines (see Figure 1-right). The intersections of these surfaces with horizontal planes are closed curves which present 6 cusp points as well as \textit{self-intersections} (see Figure 1-right and Figure 2), the arrangement of which fully characterizes the corresponding singularity of the semi-caustic. The classification of these singularities, i.e. of the distribution of the corresponding cuspidal points and self-intersections, can be synthetically represented by means of \textit{symbols}. Here, a symbol is a six-tuple of rational numbers $(s_1,...,s_6)$ where each $s_i$ is half the number of self-intersections appearing along the piece of curve joining two consecutive cusp points. Using this notation, we can state the main result of this article.

\begin{thm}[Main result]
\label{thm:MainEng}
Let $M$ be a 3-dimensional smooth and connected manifold and $\text{SubR}(M)$ be the space of contact sub-Riemannian distributions over $M$ endowed with the Withney topology. There exists an open and dense subset $\Ecal \subset \text{SubR}(M)$ such that for any $(\Delta,\g) \in \Ecal$ the following holds.
\begin{enumerate}
\item[(i)] There exists a smooth curve $\Ccal \subset M$ such that, outside $\Ccal$, the intersections of the caustic with horizontal planes $\{ h = \pm \epsilon\}$ are closed curves exhibiting 4 cusp points (see Figure 1-left below).
\item[(ii)] There exists an open and dense subset $\Ocal \subset \Ccal$ on which the intersections of the caustic with horizontal planes $\{ h = \pm \epsilon \}$ are described by pairs of symbols $(\Scal_i,\Scal_j)$ with $i,j \in \{1,2,3\}$ (see Figure 2 first and second drawings below) and
\begin{equation*}
\Scal_1 = (0,1,1,1,1,1),~ \Scal_2 = (2,1,1,1,1,1),~ \Scal_3 = (2,1,1,2,1,0).
\end{equation*}
\item[(iii)] There exists a discrete subset $\Dcal \subset \Ccal$ complement of $\Ocal$ in $\Ccal$ on which the intersections of the caustic with horizontal planes $\{ h = \pm \epsilon \}$ are described by pairs of symbols $(\Scal_i,\Scal_j)$ with $i \in \{1,2,3\}, j \in \{4,5,6,7\}$ (see Figure 2 third and fourth drawings below) and 
\begin{equation*}
\Scal_4 = ( \tfrac{1}{2},\tfrac{1}{2},1,0,0,1 ),~ \Scal_5 = (1,\tfrac{1}{2},\tfrac{1}{2},1,1,1),~ \Scal_6 = (\tfrac{3}{2},\tfrac{1}{2},1,1,0,1),~ \Scal_7 = (2,\tfrac{1}{2},\tfrac{1}{2},2,0,0).
\end{equation*}
\end{enumerate}
\end{thm}

\begin{figure}[!h]
\label{fig:1}
\hspace{-1cm}
\begin{minipage}[c]{0.33\textwidth}
\resizebox{!}{0.25 \textheight}{\begin{tikzpicture}
\node[inner sep=0pt] (A) at (0,0) {\includegraphics[scale = 0.3]{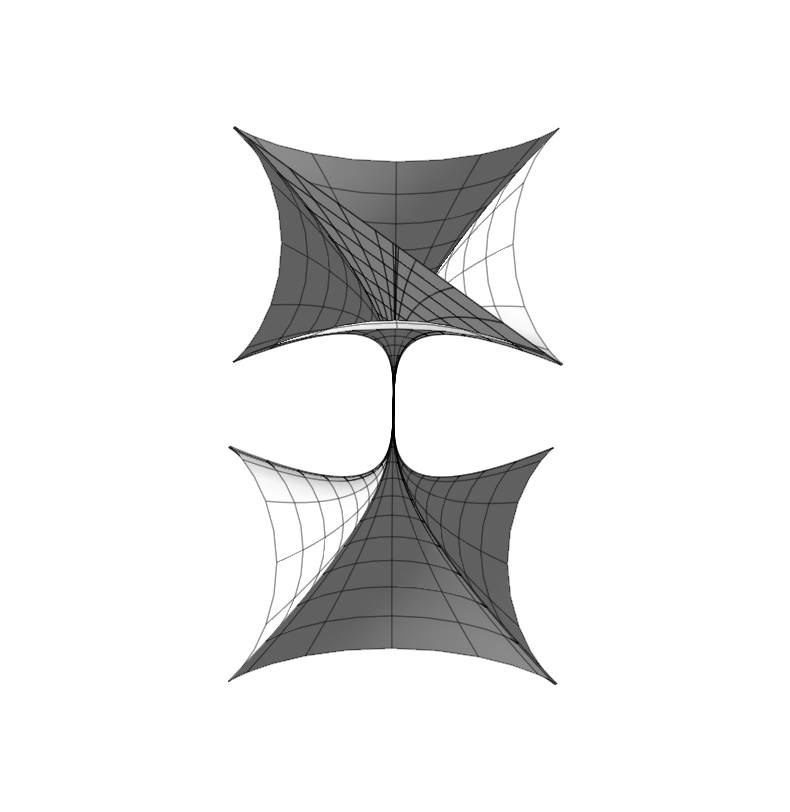}};
\node[below] at (2.5,1.6) {\footnotesize Cut Locus};
\draw[thick,->](2.5,1.6) to [out = 135 , in = 80](1,1.1);
\node[below] at (2.15,3.4) {\footnotesize Conjugate Locus};
\draw[thick,->](2.25,3.4) to [out = 135 , in = 80](0.5,2.8);
\end{tikzpicture}}
\end{minipage} 
\begin{minipage}[c]{0.33 \textwidth}
\vspace{0.8cm}
\resizebox{!}{0.25 \textheight}{\begin{tikzpicture}
\node[inner sep=0pt] (B) at (0,0) {\includegraphics[scale = 0.3]{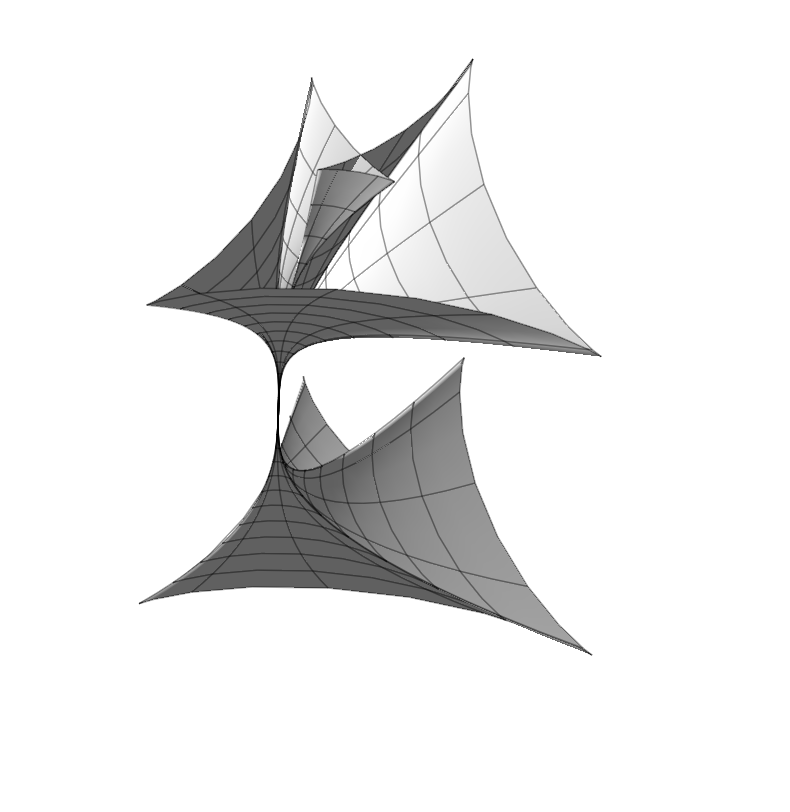}};
\node[below] at (1.75,3.75) {\footnotesize Swallow tail};
\draw[thick,->](1.75,3.75) to [out = 130 , in = 80](-0.3,2.75);
\end{tikzpicture}}
\end{minipage}
\begin{minipage}[c]{0.33 \textwidth}
\vspace{0.8cm}
\resizebox{!}{0.25 \textheight}{\begin{tikzpicture}
\node[inner sep=0pt] (B) at (0,0) {\includegraphics[scale = 0.3]{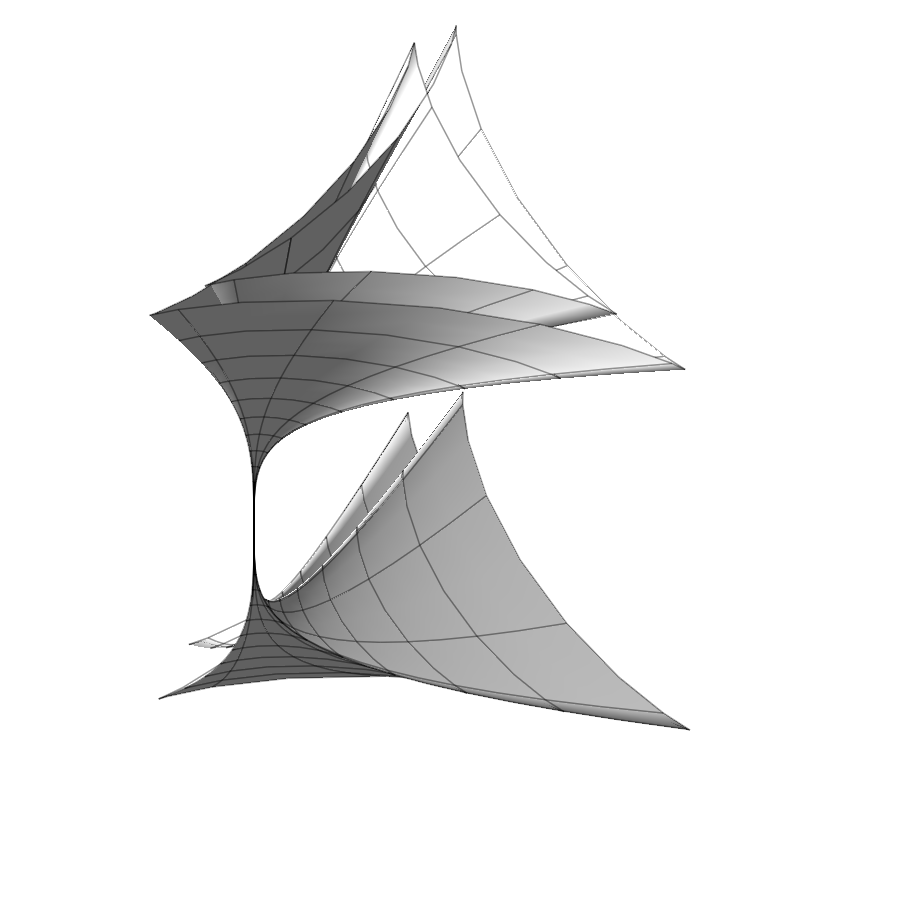}};
\draw[thick,->](3,3.75) to [out = 110 , in = 45](-0.15,3.85);
\draw[thick,->](3,3.15) to [out = 260 , in = 20](1.8,1.6);
\node[below] at (3,3.7) {\small Self-intersections};
\draw[thick,->]  (1.25,0.25) to [out = 180 , in = 270] (0.75,0.75);
\node[right] at (1.25,0.25) {\small Plea line};
\end{tikzpicture}}
\end{minipage}
\vspace{-0.75cm} \caption{Generic 4-cusp conjugate and cut loci outside $\Ccal$ (left), formation of the swallow tail near $\Ccal$ (center), folding of the swallow tail and formation of the 6-cusp singularity along $\Ccal$ (right)}
\end{figure}
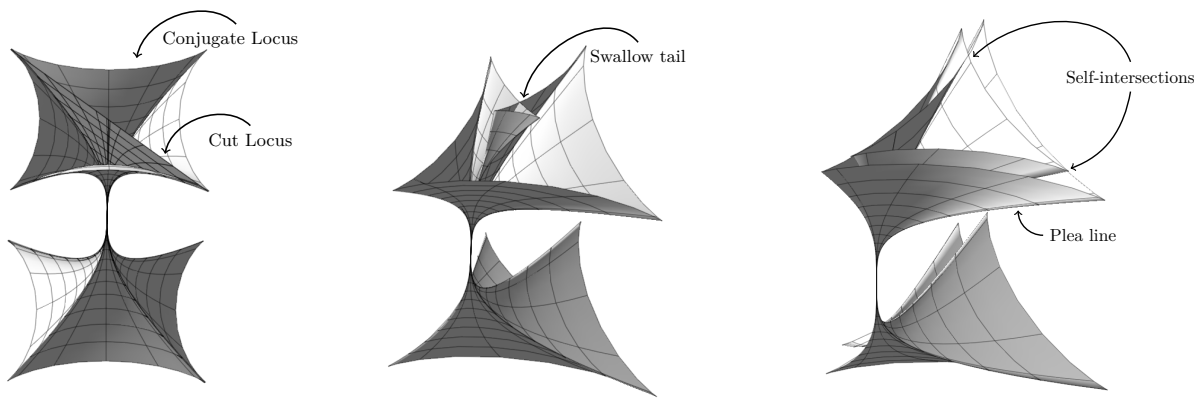

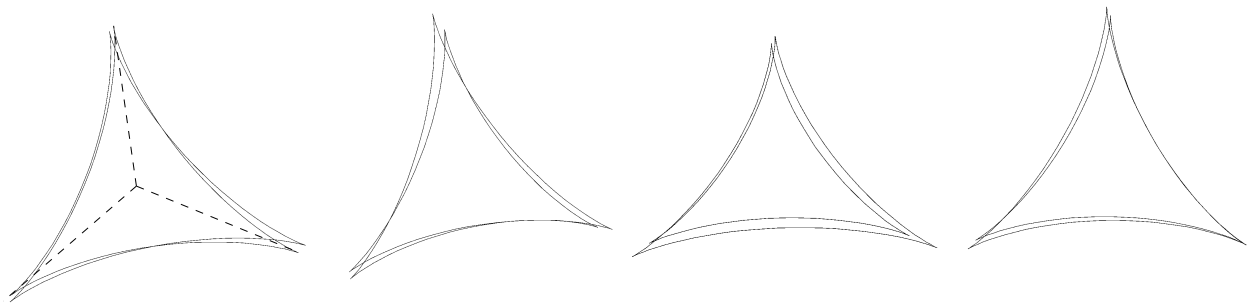
\begin{figure}[!h]
\label{fig:3}
\centering      
      \begin{tikzpicture}
      \node at (-3.5,0) {\includegraphics[scale = 0.11]{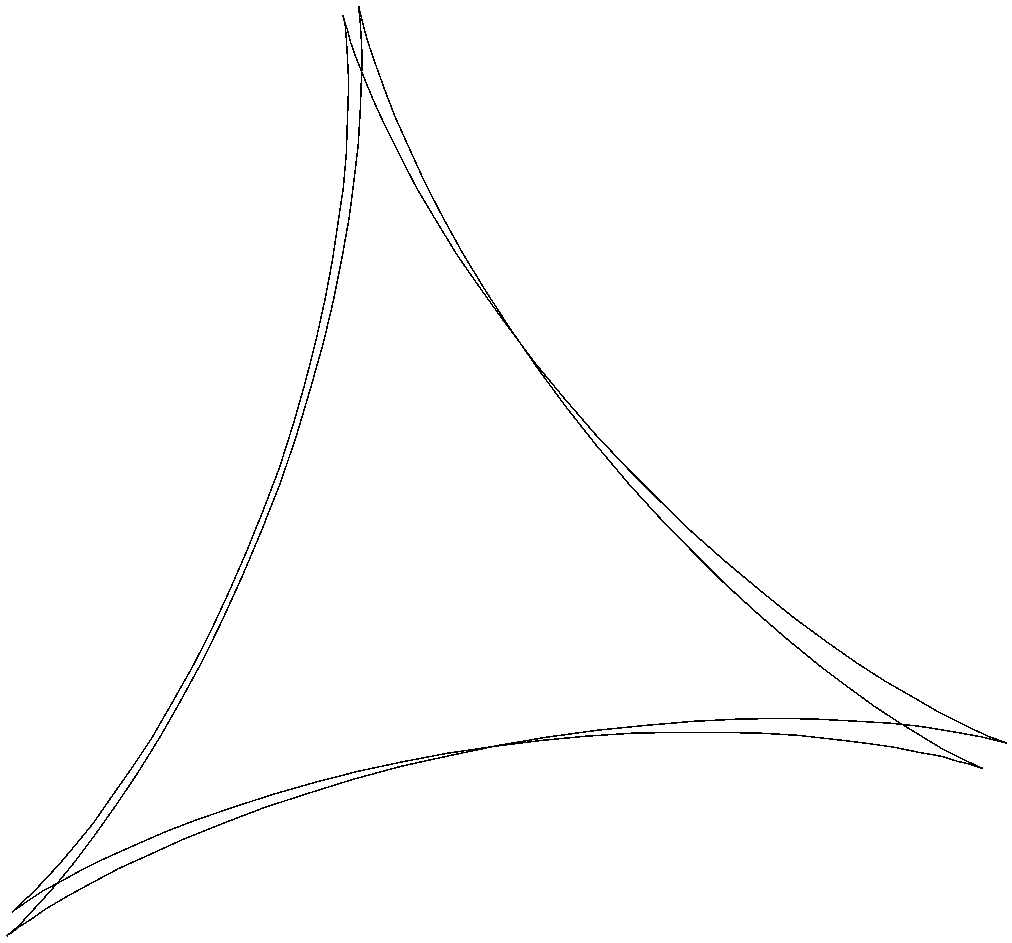}};
      \draw[very thin,dashed](-3.775,-0.275)--(-4.075,1.85);
      \draw[very thin,dashed](-3.775,-0.275)--(-1.65,-1.15);
      \draw[very thin,dashed](-3.775,-0.275)--(-5.525,-1.8);
      \node at (0.75,0.25) {\includegraphics[scale = 0.11]{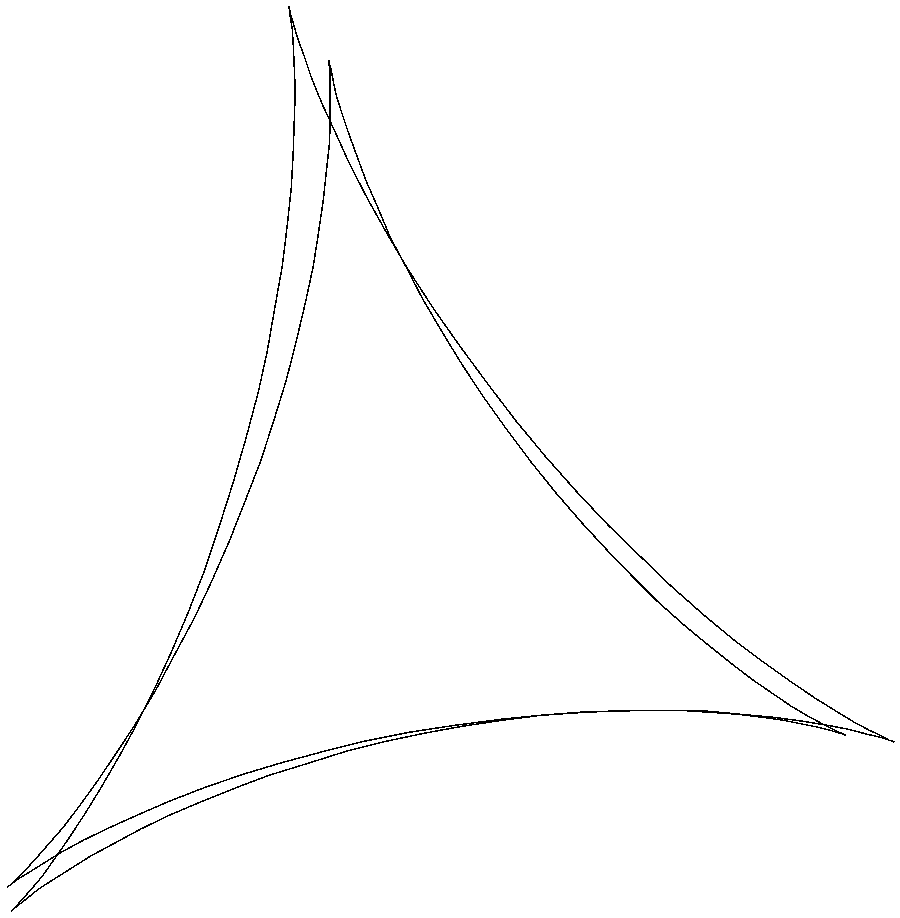}};
      \node at (4.75,0.25) {\includegraphics[scale = 0.125]{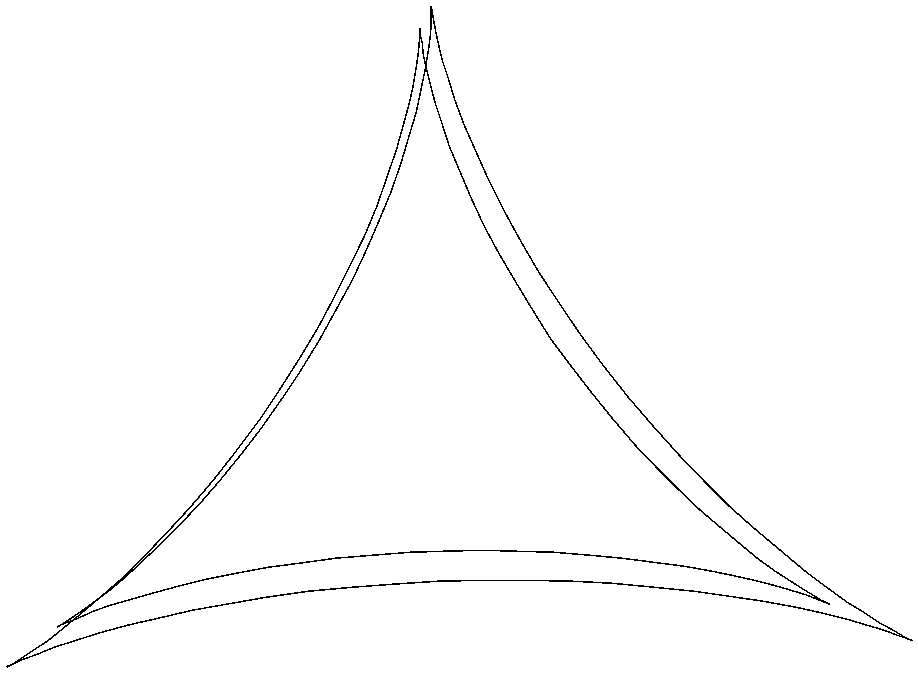}};
      \node at (9,0.5) {\includegraphics[scale = 0.115]{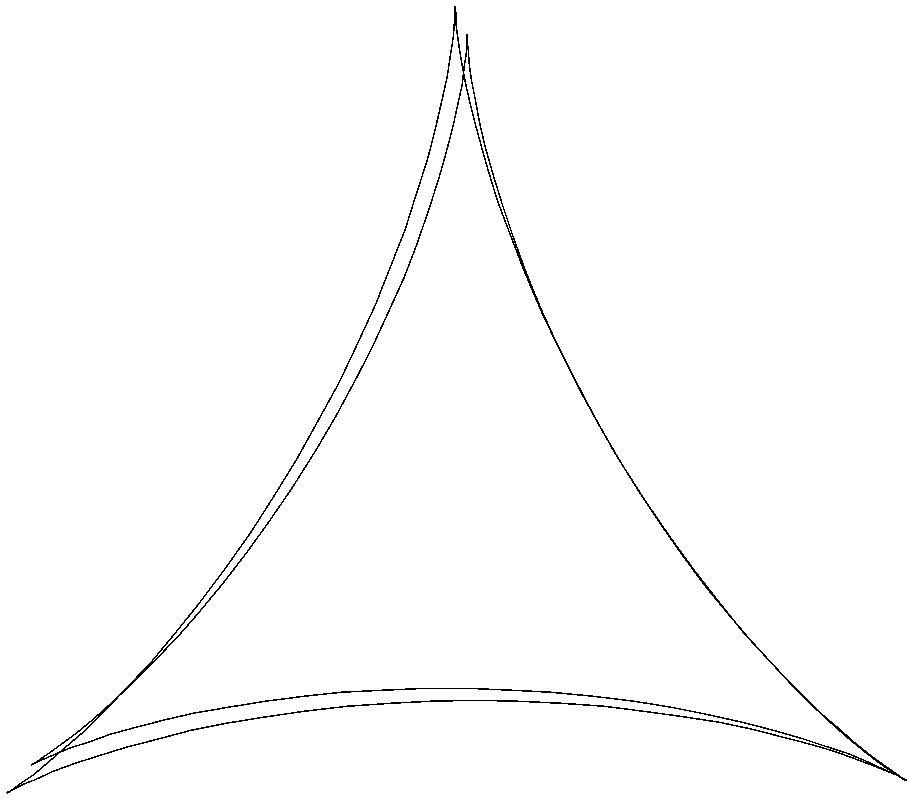}};
	  \end{tikzpicture}
\caption{Intersections of the upper semi-caustics corresponding to the symbols $\Scal_1,\Scal_2,\Scal_4$ and $\Scal_6$ with the planes $h = \epsilon$, and intersection of the corresponding cut locus (dashed lines on the left)}
\end{figure}

This result is in a sense the most natural one to be expected after the classification of the half conjugate loci displayed in \cite{ACG}. Indeed, it transcribes the fact that the upper and lower semi-caustics are independent and that there are no extra couplings appearing between the two structures. Indeed, the only possible obstruction to the combination of two given symbols is that the corresponding codimension in the space of Taylor coefficients -- which is preserved by standard arguments of transversality theory -- is strictly larger than 3. In particular, this generically prevents pairs of the form $(\Scal_i,\Scal_j)$ with $i,j \in \{4,5,6,7\}$ from appearing. 

\section{3D-Contact sub-Riemannian manifolds and their conjugate locus}

In this section, we recall some elementary facts about sub-Riemannian geometry defined over $3$-dimensional manifolds. For a complete introduction, see e.g. \cite{ABB}.

\begin{Def}[Sub-Riemannian manifold]
\label{def:3DContact}
A 3D-contact sub-Riemannian manifold is defined by a triple $(M,\Delta,\g)$ where 
\begin{enumerate}
\item[$\bullet$] $M$ is an $3$-dimensional smooth and connected differentiable manifold,
\item[$\bullet$] $\Delta$ is a smooth $2$-dimensional distribution over $M$ with step 1, i.e.  
\begin{equation*}
\Span \left\{ X_1(q),X_2(q),[X_1(q),X_2(q)] \right\} = T_q M, 
\end{equation*}
for all $q \in M$ and $(X_1(q),X_2(q))$ spanning $\Delta(q)$. 
\item[$\bullet$] $\g$ is a Riemannian metric over $\Delta$. 
\end{enumerate}
\end{Def}

\begin{Def}[Horizontal curves and sub-Riemannian metric]
An absolutely continuous curve $\gamma(\cdot)$ is said to be \textnormal{horizontal} if $\dot \gamma(t) \in \Delta(\gamma(t))$ for $\Lcal^1$-almost every $t \in [0,T]$. We define the \textnormal{length} $l(\gamma(\cdot))$ of a horizontal curve $\gamma(\cdot)$ as
\begin{equation*}
l(\gamma(\cdot)) = \INTSeg{\sqrt{\g_{\gamma(t)}(\dot \gamma(t),\dot \gamma(t)) \,}}{t}{0}{T}.
\end{equation*}
For $(q_0,q_1) \in M$, it is then possible to define the \textnormal{sub-Riemannian distance} $d_{\textnormal{sR}}(q_0,q_1)$ as the infimum of the length of the horizontal curves connecting $q_0$ and $q_1$. 
\end{Def}

Given a local orthonormal frames $(X_1,X_2)$ for the metric $\g$ which spans $\Delta$, the Carnot-Carath\'eodory distance $d_{\textnormal{sR}}(q_0,q_1)$ can be alternatively computed by solving the optimal control problem
\begin{equation}
\label{eq:OCP}
\left\{
\begin{aligned}
& \hspace{-0.7cm} \min_{(u_1(\cdot),u_2(\cdot))} \INTSeg{(u_1^2(t)+u_2^2(t))}{t}{0}{T} \\
\text{s.t.} & ~~ \dot \gamma(t) = u_1(t) X_1(\gamma(t)) + u_2(t)X_2(\gamma(t)), \hspace{0.2cm}  u_1^2(t) + u_2^2(t) \leq 1, \\
\text{and} & ~~ (\gamma(0),\gamma(T)) = (q_0,q_1).
\end{aligned}
\right.
\end{equation}
We detail in the following proposition the explicit form of 3D-contact sub-Riemannian geodesics obtained by applying the maximum principle to \eqref{eq:OCP}.

\begin{prop}[The Pontryagin Maximum Principle in the 3D contact case] 
\label{prop:PMP_Contact} \hfill \\
Let $\gamma(\cdot) \in \Lip([0,T],M)$ be a horizontal curve and $\Hpazo : T^*M \rightarrow \R$ be the Hamiltonian associated to the contact geodesic problem, defined by 
\begin{equation*}
\Hpazo (q,\lambda) = \tfrac{1}{2} \langle \lambda , X_1(q) + X_2(q) \rangle,
\end{equation*}
for any $(q,\lambda) \in T^*M$. Then, the curve $\gamma(\cdot)$ is a contact geodesic parametrized by sub-Riemannian arclength if and only if there exists a Lipschitzian curve $t \in [0,T] \mapsto \lambda(t) \in T^*_{\gamma(t)} M$ such that $t \mapsto (\gamma(t),\lambda(t))$ is a solution of the Hamiltonian system
\begin{equation}
\label{eq:PMP_Contact}
\dot \gamma(t) = \partial_{\lambda}\Hpazo(\gamma(t),\lambda(t)), \qquad \dot \lambda(t) = -\partial_q \Hpazo(\gamma(t),\lambda(t)), \qquad \Hpazo(\gamma(t),\lambda(t)) = \tfrac{1}{2}.
\end{equation}
We denote by $\overrightarrow{\Hpazo} \in \Vecf(T^* M)$ the corresponding Hamiltonian vector field defined over the cotangent bundle and by $(\gamma(t),\lambda(t)) = e^{t \overrightarrow{\Hpazo}}(q_0,\lambda_0)$ the corresponding solution of \eqref{eq:PMP_Contact}.
\end{prop}

Both the absence of abnormal lifts and the sufficiency of the maximum principle are consequences of the contact hypothesis made on the sub-Riemannian structure (see e.g. \cite[Chapter 4]{ABB}). 

\begin{Def}[Exponential map]
Let $q_0 \in M$ and $\Lambda_{q_0} = \{ \lambda \in T^*_{q_0} M ~\text{s.t.}~ \Hpazo(q_0,\lambda_0) = \frac{1}{2} \}$. We define the \textnormal{exponential map} from $q_0$ as 
\begin{equation*}
\E_{q_0} : (t,\lambda_0) \in \R_+ \times \Lambda_{q_0} \mapsto \pi_M \left( e^{t \overrightarrow{\Hpazo}}(q_0,\lambda_0) \right),
\end{equation*}
where $\pi_M : T^*M \rightarrow M$ is the canonical projection. 
\end{Def}

In this paper, we study the germ at the origin -- i.e. equivalence classes of maps defined by equality of derivatives up to a certain order -- of the \textit{cut} and \textit{conjugate loci} associated to contact sub-Riemannian structures.

\begin{Def}[Cut and conjugate locus]
\label{def:CL}
Let $(M,\Delta,\g)$ be a 3D contact sub-Riemannian manifold, $(q_0,\lambda_0) \in T^*M$  and $\gamma(\cdot) = \E_{q_0}(\cdot,\lambda_0)$ be a geodesic parametrized by sub-Riemannian arclength. The \textnormal{cut time} associated to $\gamma(\cdot)$ is defined by 
\begin{equation*}
\tau_{\textnormal{cut}} = \sup \{ t \in \R_+ ~\text{s.t. $\gamma_{[0,t)}(\cdot)$ is optimal } \},
\end{equation*}
and the corresponding \textnormal{cut locus} is 
\begin{equation*}
\textnormal{Cut}(q_0) = \{ \gamma(\tau_{\textnormal{cut}}) ~\text{s.t. $\gamma(\cdot)$ is a sub-Riemannian geodesic from $q_0$}\}.
\end{equation*}
The \textnormal{first conjugate time} $\tau_{\textnormal{conj}}$ associated to the curve $\gamma(\cdot)$ is define by 
\begin{equation*}
\tau_{\textnormal{conj}} = \inf \{ t \in \R_+ ~\text{s.t. $(t,\lambda_0)$ is a critical point of $\E_{q_0}(\cdot,\cdot)$} \},
\end{equation*}
and the corresponding \textnormal{conjugate locus} is then defined by
\begin{equation*}
\Conj(q_0) = \{ \gamma(\tau_{\textnormal{conj}}) ~\text{s.t. $\gamma(\cdot)$ is a sub-Riemannian geodesic from $q_0$}\}.
\end{equation*}
\end{Def}

We recall in Theorem \ref{thm:NormalForm} below the \textit{normal form} introduced formally in \cite{CGK} and then derived geometrically in \cite{ACG} for 3D-contact sub-Riemannian structures. Up to a simple change of coordinates, we can assume that $0 \in M$ and study the germ of the conjugate locus in a neighbourhood of the origin.

\begin{thm}[Normal form for 3D-contact sub-Riemannian distributions]
\label{thm:NormalForm}
Let $(M,\Delta,\g)$ be a 3D-contact sub-Riemannian structure and $(X_1,X_2)$ be a local orthonormal frame for $\Delta$ in a neighbourhood of the origin. Then, there exists a smooth system of so-called \textnormal{normal coordinates} $(x,y,w)$ on $M$ along with two maps $\beta,\gamma \in C^{\infty}(M,\R)$ such that $(X_1,X_2)$ can be written in \textnormal{normal form} as
\begin{equation}
\left\{
\begin{aligned}
X_1(x,y,w) & = (1+y^2 \beta(x,y,w)) \partial_x - xy \beta(x,y,w) \partial_y + \tfrac{y}{2} (1+\gamma(x,y,w)) \partial_w, \\
X_2(x,y,w) & = (1+x^2 \beta(x,y,w)) \partial_y - xy \beta(x,y,w) \partial_x - \tfrac{x}{2} (1+\gamma(x,y,w)) \partial_w, \\
\beta(0,0,w) & = \gamma(0,0,w) = \partial_x \gamma(0,0,w) = \partial_y \gamma(0,0,w) = 0.
\end{aligned}
\right.
\end{equation}
This system of coordinates is unique up to an action of $SO(2)$ and adapted to the contact structure, i.e. it induces a gradation with respective weights  $(1,1,2)$ on the space of formal power series in $(x,y,w)$. 
\end{thm}

The truncation at the order $(-1,-1)$ of this normal form is precisely the orthonormal frame associated with the usual left-invariant metric on the Heisenberg group, i.e. $(X_1,X_2) = (\partial_x +\frac{y}{2}\partial_w,\partial_y - \frac{x}{2} \partial_w)$. 

Given a Heisenberg geodesic with initial covector $\lambda_0 = (p(0),q(0),r(0)) \in T^*_0M$, the corresponding conjugate time is exactly $\tau_{\textnormal{conj}} = 2 \pi /r(0)$. In the general contact case (see e.g. \cite[Chapter 16]{ABB}), the first conjugate time is of the form 
\begin{equation*}
\tau_{\textnormal{conj}} = 2 \pi /r(0) + O(1/r(0)^3)
\end{equation*}
where the higher order terms can be expressed via the coefficients of the Taylor expansions
\begin{equation*}
\beta(x,y,w) = \sum\nolimits_{l=1}^k \beta^l(x,y,w) + O^{k+1}(x,y,w), \qquad \gamma(x,y,w) = \sum\nolimits_{l=2}^k \gamma^l(x,y,w) + O^{k+1}(x,y,w)
\end{equation*}
of the maps $\beta$ and $\gamma$ with respect to the gradation $(1,1,2)$. We introduce in the following equations the coefficients $(c_i,c_{jk},c_{lmn})$ of the polynomial functions $\gamma^2,\gamma^3$ and $\gamma^4$ appearing in these expansions. 
\begin{equation}
\label{eq:NotaTaylor}
\left\{
\begin{aligned}
\gamma^2(x,y) & = (c_0 + c_2)(x^2 + y^2) + (c_0 - c_2)(x^2 - y^2) - 2 c_1 xy ~,~ \\
\gamma^3(x,y) & = (c_{11}x + c_{12}y)(x^2+y^2) + 3(c_{31}x - c_{32}y)(x^2-y^2) - 2(c_{31}x^3+c_{32}y^3) \\
\gamma^4(x,y,w) & = \frac{w}{2} \Big( (c_{421} + c_{422})(x^2 + y^2) + (c_{421} - c_{422}) (x^2 - y^2) - 2c_{423}xy \Big)  + c_{441} (x^2 + y^2)^2 \\
& + c_{442} (x^4 + y^4 - 6x^2 y^2) + 4 c_{443} \, xy (x^2 - y^2) + c_{444} (x^4 - y^4) - 2 c_{445} \, xy(x^2 + y^2)
\end{aligned}
\right.
\end{equation}
These coefficients derive from the irreducible decompositions of $\gamma^2$, $\gamma^3$ and $\gamma^4$ under the action of $SO(2)$, and are precisely the fundamental invariants discriminating the singularities listed in Theorem \ref{thm:MainEng}. Indeed, it was proven in \cite{CGK} that the caustic becomes degenerate along a smooth curve $\Ccal \subset M$ on which the coefficients $(c_0,c_1,c_2)$ in the decomposition of $\gamma^2$ in \eqref{eq:NotaTaylor} satisfy $c_0 = c_2$ and $c_1 = 0$. Moreover, the decompositions of $\beta,\gamma^2$ and $\gamma^3$ have an interpretation in terms of canonical sub-Riemannian \textit{connection} and \textit{curvature}, as detailed in \cite{CGK}. 

Following the methodology developed in \cite{ACG}, we introduce the coordinates $(h,\varphi)$ defined by
\begin{equation*}
(p,q) = (\cos(\varphi),\sin(\varphi)),~ h = \textnormal{sign}(w) \sqrt{|w|/\pi}.
\end{equation*}
We can then express the Taylor expansion of order $k \geq 3$ of the semi conjugate loci as the suspension
\begin{equation}
\label{eq:Suspension}
\Conj_{\pm}(\varphi,h) = (x(\varphi,h),y(\varphi,h),h) = \left( \sum\nolimits_{l=3}^k h^l f_l^{\pm}(\varphi),h \right),
\end{equation}
where the $\pm$ symbol highlights the dependency of the expression on the sign of $w$. By carrying out explicitly the computations necessary to put the conjugate locus in the suspended form \eqref{eq:Suspension}, it can be verified that $f_3^+ = f_3^-$ and $f_4^+ = f_4^-$. Hence, the generic behaviour of the full conjugate locus is already known outside the smooth curve $\Ccal \subset M$ on which the caustic becomes degenerate. We therefore restrict our attention to the more degenerate singularities arising in the form of \textit{self-intersections} along the curve $\Ccal$.

To prove that the symbols listed in Theorem \ref{thm:MainEng} are generic, one needs to compute the Taylor expansions of order $k=7$ of the metric in a neighbourhood of the origin. All these expressions where obtained using Maple software using a piece of code that can be found at the following address:

\begin{center}
\verb?http://www.lsis.org/bonnetb/depots/Sub-Riemannian_Conjugate_Locus?
\end{center}

\vspace{-0.25cm}


\section{Generic singularities of the full-conjugate locus}

In this section, we prove Theorem \ref{thm:MainEng}. We show that the coefficients defined in \eqref{eq:NotaTaylor}, characterizing the generic conjugate locus on the curve $\Ccal$, generate independent structures for $\Conj_+$ and $\Conj_-$. 

\begin{Def}[Self-intersection set of the semi-conjugate loci]
We define the self-intersection set of $\Conj_+$ as
\begin{equation*}
\Self(\Conj_+) = \{ (h,\varphi_1,\varphi_2) \in \R_+ \times \Sbb^1 \times \Sbb^1 ~\text{s.t.}~ h > 0,~ \varphi_1 \neq \varphi_2,~ \Conj_+(h,\varphi_1) = \Conj_+(h,\varphi_2) \}.
\end{equation*}
An angle $\varphi \in \Sbb^1$ is said to be \textnormal{adherent} to $\Self(\Conj_+)$ provided that $(0,\varphi,\varphi+\pi) \in \overline{\Self(\Conj_+)}$. The set of such angles is denoted by $\ASelf(\Conj_+)$. We define in the same way the self-intersection and adherent angles sets of $\Conj_-$. 
\end{Def}

We recall in Theorem \ref{thm:ACG} below the complete classification of the generic self-intersections of $\Conj_+$ (or equivalently of $\Conj_-$) which was derived in \cite{ACG}. 

\begin{thm}[Generic self-intersections of the positive semi conjugate locus]
\label{thm:ACG}
Let $\Ccal \subset M$ be the curve defined in Theorem \ref{thm:MainEng}. Then, $\Ccal$ is a generically smooth curve, and outside $\Ccal$ the semi conjugate loci are the standard 4-cusp semi-caustics. On $\Ccal$, the following situations can occur. 
\begin{enumerate}
\item[(i)] There exists an open and dense subset $\Ocal^+ \subset \Ccal$ on which the self-intersections of $\Conj_+$ are described by the symbols $(\Scal_1,\Scal_2,\Scal_3)$.
\item[(ii)] There exists a discrete subset $\Dcal^+$ complement of $\Ocal^+$ in $\Ccal$  on which the self-intersections of $\Conj_+$ are described by the symbols $(\Scal_4,\Scal_5,\Scal_6,\Scal_7)$.
\end{enumerate}
\end{thm}

We recall in the following lemma the structural result allowing to describe the sets $\ASelf(\Conj_{\pm})$.

\vspace{0.1cm}
\begin{lem}[Structure of the self-intersection]
The set of adherent angle to $\Self(\Conj_{\pm}$ satisfies the inclusion $\ASelf(\Conj_{\pm}) \subset \{ \varphi \in \Sbb^1 ~\text{s.t.}~ f'_4(\varphi) \wedge f_5^{\pm}(\varphi) = 0 \}$, where $a \wedge b \equiv \textnormal{det}(a,b)$ for vectors $a,b \in \R^2$. 
\end{lem}
\vspace{0.1cm}

An explicit computation based on the expressions of the maps $f_4,f_5^-$ and $f_5^+$ defined in \eqref{eq:Suspension} shows that 
\begin{equation*}
f'_4(\varphi) \wedge f_5^{\pm}(\varphi) = -20\pi^2 \tilde{b} \sin(3\varphi + \omega_b) P_{\pm}(\varphi)
\end{equation*}
where we introduced the notations $b = (c_{31},c_{32}) = \tilde{b}(\sin(\omega_b),-\cos(\omega_b)) \in \C$ and 
\begin{equation*}
P_{\pm}(\varphi) = A_{\pm} \cos(2\varphi) + B_{\pm} \sin(2\varphi) + C \cos(4 \varphi) + D \sin(4 \varphi).
\end{equation*}
Here, the coefficients $(A_{\pm},B_{\pm},C,D)$ are \textbf{independent linear combinations} of the fourth-order coefficients $(c_{421},c_{422},c_{423},c_{442},c_{443},c_{444},c_{445})$ introduced in \eqref{eq:NotaTaylor}. Their analytical expressions were derived again by using the Maple software and write as follows:
\begin{equation}
\label{eq:CoeffExpression}
\left\{
\begin{aligned}
& A_{\pm} = \tfrac{35}{8}(c_{422} - c_{421}) \pm 3 \pi c_{423} + 45 c_{445},~  C = 36 c_{443}, \\
& B_{\pm} = \tfrac{35}{8} c_{423} \pm 3 \pi (c_{421} - c_{422}) + 45 c_{444},~ D = -36 c_{442}.
\end{aligned}
\right.
\end{equation}

In order to describe the roots of $f'_4 \wedge f_5^{\pm}$ in $\Sbb^1$, it is convenient to introduce the complex polynomials 
\begin{equation*}
\tilde{P}_{\pm}(z) = \mu z^4 + \nu_{\pm} z^3 + \bar{\nu}_{\pm} z + \bar{\mu},~~ \tilde{T}(z) = b z^3 + \bar{b},
\end{equation*}
where the complex coefficients $\nu_-$, $\nu_+$ and $\mu$ are defined by 
\begin{equation*}
\nu_{\pm} = \tfrac{1}{2} (A_{\pm} - i B_{\pm}), ~~ \mu = \tfrac{1}{2} (C - iD).
\end{equation*}
It can be checked that $f_4'(\varphi) \wedge f_5^{\pm}(\varphi) = 0$ if and only if $\tilde{P}_{\pm}(e^{i \varphi}) \tilde{T}(e^{i \varphi}) = 0$. Since $(A_{\pm},B_{\pm},C,D)$ independently span $\R^6$ as a consequence of \eqref{eq:CoeffExpression}, the classification of the generic singularities of the degenerate semi-caustics reduces to understanding the distribution of the unit roots of $\tilde{P}_{\pm}$ and $\tilde{T}$ as $(\mu,\nu_{\pm},b)$ span $\C^4$. In the following lemma, we compile some of the facts highlighting the relationship between these roots and the symbols describing $\Self(\Conj_{\pm})$ that can be found in \cite{ACG,CGK}.

\begin{lem}[Algebraic equations describing the singularities]
Let $\gamma^2,\gamma^3$ and $\gamma^4$ be given as in \eqref{eq:NotaTaylor} along with their decomposition under the action of $SO(2)$. Then, the following hold.
\begin{enumerate}
\item[(i)] The polynomials $\tilde{P}_{\pm}$ have either 2 or 4 distinct simple roots on the unit circle.
\item[(ii)] The polynomials $\tilde{P}_{\pm}$ share a unit root with $\tilde{T}$ if and only if  $\textnormal{Res}(\mu,\nu_{\pm},b) = 0$ where $\textnormal{Res}(\mu,\nu_{\pm},b)$ is the resultant polynomial of $\tilde{P}_{\pm}$ and $\tilde{T}$. 
\end{enumerate}
\end{lem} 

The symbols $\Scal_i$ introduced in Section \ref{section:Introduction} can be understood as follows: $\Scal_1,\Scal_2,\Scal_3$ describe all the possible situations in which $\tilde{P}_{\pm}$ have 2 or 4 simple roots, none of which is shared with $\tilde{T}$, while $\Scal_4,\Scal_5,\Scal_6,\Scal_7$ refer to all the combinations in which $\tilde{P}_{\pm}$ and $\tilde{T}$ share a single simple root. All the other possible situations are non-generic even for the semi conjugate loci.

For the full conjugate locus, the situation is the following. It can be shown that the sets $S_{\pm} =  \{(c_{jk},c_{lmn}) \in \R^{16} ~\text{s.t.}~ \mu,b \neq 0,~ \textnormal{Res}(\mu,\nu_{\pm},b) = 0 \}$ are smooth manifolds of codimension 1. Moreover, it holds that $S_- \cap S_+$ has codimension 2 as a consequence of the fact that $\textnormal{Res}(\mu,\nu_{\pm},b)$ is a homogeneous polynomial of degree 7 in $(\mu,\nu_{\pm},b)$ which only vanishes at isolated points where its differential is surjective when both $b,\mu  \neq 0$.

Since $M$ is a 3-dimensional manifold and $\Ccal \subset M$ already has codimension 2, it can be proven that the set of sub-Riemannian distributions with canonical orthonormal frame $(X_1,X_2)$ defined in \eqref{eq:NotaTaylor} such that $(c_{jk},c_{lmn}) \notin S_- \cap S_+$ is open and dense. This stems from standard transversality and preservation of codimension arguments, following \cite[Section 3.3]{ACG}. Hence, it cannot happen generically that the singularities of the conjugate locus are described by a pair of symbols $(\Scal_i,\Scal_j)$ with both $i,j \in \{4,5,6,7\}$.

This independence result implies in particular the following structural corollary for the cut locus. 

\begin{cor}[Generic behaviour of the 3D-contact sub-Riemannian cut locus]
\label{cor:Cut}
There exists an open and dense subset $\Ecal \subset \text{SubR}(M)$ for the Whitney topology such that for any $(\Delta,\g) \in \Ecal$, the following holds. 

\begin{enumerate}
\item[(i)] Outside the smooth curve $\Ccal \subset M$ defined in Theorem \ref{thm:MainEng}, the semi cut loci are independent : each of the two is a portion of plane joining two opposite plea lines along the corresponding semi caustic (see Figure 1-left). 
\item[(ii)] On the curve $\Ccal$, the semi cut loci are the union of three portions of planes connecting the $h$-axis with alternate plea lines along the corresponding semi caustic (see Figure 2-left dotted lines). There is no interdependence whatsoever between these planes for the positive and negative semi cut loci.
\end{enumerate}
\end{cor}

The generic behaviour of the 3D-contact semi cut loci outside the degenerate curve $\Ccal$ is already known (see e.g. \cite[Chapter 16]{ABB}). As described in the introduction, swallow tails appear along the semi caustics as one approaches the curve $\Ccal$ (see Figure 1-center), the four plea lines then degenerating into the 6-cusp structure. The corresponding cusps are regrouped by pairs of the form $(k,k+1)$, appearing respectively in a small vicinity of the cuspidal angles $\tfrac{k}{3}(\pi-\omega_b)$ with $k \in \{0,1,2\}$. 

Generically, the semi cut loci are the unions of three portions of planes joining three of the six plea lines to the $h$-axis (see Figure 2-left dotted lines). Moreover, in the absence of an interior cusp, these supporting plea lines can be chosen freely for each semi cut loci to be either one of the two numberings $(1,3,5)$ or $(2,4,6)$ of the six plea lines. The result of this paper stating that the cusps of the conjugate locus are independently distributed for its upper and lower parts therefore yields Corollary \ref{cor:Cut}. 

\section*{Acknowledgements}

\noindent \scriptsize{\textit{This research is partially supported by the University of Padova grant SID 
2018 ``Controllability, stabilizability and infimun gaps for control 
systems'', prot. BIRD 187147, by the Archim\`ede
Labex (ANR-11-LABX-0033), by the A*MIDEX project (ANR-
11-IDEX-0001-02) and by the SRGI ANR Grant ANR-15-CE40-0018}}

\end{document}